\documentclass[a4paper,12pt]{amsart} 
\usepackage{amscd}             
\usepackage{amssymb}
\usepackage{amsthm}
\usepackage{epsf} 
\usepackage[T1]{fontenc}            
\newtheorem{thm}{Theorem}[section]

\newtheorem{prop}[thm]{Proposition}
\newtheorem{example}[thm]{Example}
\newtheorem{defn}[thm]{Definition}

\newtheorem{rem}[thm]{Remark}

\def\c1{\operatorname{c_1}}
\def\c2{\operatorname{c_2}}

\def\PP{{\mathbf P}}

\def\O{{\mathcal O}}

\def\H{{\mathcal H}}

\def\+{\oplus}                   % direct sum
\def\*{\otimes}                  % tensor product
\begin{document}
\thispagestyle{empty}
 \vspace*{3.5 cm}
 \centerline{\bf RANK TWO BUNDLES ON ALGEBRAIC CURVES}
   \centerline{\bf AND DECODING OF GOPPA CODES}
% \title{Rank two bundles on Algebraic Curves and Decoding of Goppa Codes}
%\maketitle
  \bigskip  \bigskip
  \bigskip

\centerline{Trygve Johnsen}  
  
  \centerline{Department of Mathematics, University of Bergen, } 
  \centerline{Johs. Bruns gt. 12, N-5008 Bergen, NORWAY}
\centerline{e-mail: johnsen@mi.uib.no}
\medskip
  \bigskip  \bigskip
  \bigskip
\noindent {\bf Abstract}
\noindent We study a connection between two topics: Decoding of Goppa codes
arising from an algebraic curve, and rank two extensions of
certain line bundles on the curve. The material about each isolated
topic is well known. Our contribution is just to expose
a connection between them. 

\noindent {\bf AMS Subj. Classification:} 14J28 (14H45).


% Section 1
\section{\bf Introduction}
\label{intro}

Let  $C$  be a  curve over a finite field $F_q$. Some years ago
V. Goppa showed how to produce codes from such a curve (For a survey,
see for example Gieseker \cite {G1} ). In this note we will show how socalled
syndrome decoding of (the duals of the original ) Goppa codes in an
intimate way is connected to the study of rank two bundles, that are
extensions of the structure sheaf $\O_C$ and a certain fixed line
bundle on $C$ determined by the code in question. In fact, each
syndrome corresponds to an extension.
    Moreover syndromes due to correctable error vectors always
correpond to unstable rank two bundles, that is: semistable bundles
never correspond to syndromes of correctable errors (Here we call
an error vector, or simply error, correctable if its Hamming weight
is at most $\frac{d-1}{2}$, where $d$ is the designed minimum distance of
the code). For correctable errors the process of error location is
translated into finding a certain quotient line bundle of minimal
degree (or dually: a line subbundle of maximal degree) of the rank two
bundle defined by the syndrome, and then pick
the relevant section of this line bundle.
If an error is not correctable, there is not necessarily a unique such 
quotient line bundle of minimal degree.
   We gain our insight  through a certain projective embedding
$C \subset  \PP $, with the property that the columns of the parity
check matrix of the Goppa code in question are interpreted as points
of the embedded copy of $C$. We study  $j$-secant $(j-1)$-planes to
$C$, for $j =1,2,...$ . Using the proper definitions of these geometrical
objects, we thus obtain a stratification of $\PP$, which viewed from
one angle is a stratification of syndromes, according to how many
errors that have to be made to obtain the syndrome. Viewed from
another angle it is a stratification of rank two extensions according
to the socalled $s$-invariant.
    The geometric picture associated to rank two extensions was
given very explicitly in Lange and Narasimhan \cite{LN}, and the basic
idea was given already in Atiyah \cite{A}. We just remark that the
aspects interesting
to us carry over to the case of positive characteristic.
Then we compare with the picture obtained from the space of syndrome
vectors.
    So far we have not been able to utilize these geometrical
observations to make any constructive decoding algorithms. To do
so one would have to introduce some kind of "extention
arithmetic" to perform decoding.  Some explicit considerations along
these lines are made in \cite{BC}. For constructive decoding
algorithms for Goppa codes in general, see the celebrated Feng and Rao \cite{FR}, or for
example Duursma \cite{D}, Justesen et al \cite{JLJH}, Pellikaan
\cite{P}, Skorobogatov and Vladut \cite{SV}, or the
special issue IEEE Trans. of Info. Theory, Vol. 41, No.6, Nov. 1995.

   In Section 2 we recall some basic facts concerning
algebraic-geometric (Goppa) codes. We also define and make an
elementary study of the secant varieties that in a natural way
turn up in connection with these codes.
In Section 3 we introduce vector bundle language, and in Theorem 
\ref{mainthm} we present the connection between correctable errors and
unstable bundle extensions of rank two.

\begin{rem} \label{CIproj}
{\rm In Goppa \cite{G2}  one describes an algorithm for decoding of Goppa
codes from rational curves. There one assumes first of all that at
most  $t$  errors are made, where $t = [(d-1)/2]$. Then one assumes
that exactly $t$ errors are made, and sets up a system of equations
to solve the problem of error-location given that this extra
assumption holds. If the equations yield no solution, one assumes
that $t-1$ errors are made, sets up a new system of equations,
and so on. In the end one arrives at a point where one finds a
solution since the basic assumption is that at most  $t$  errors
are made. The coefficient matrices, set up to find the elementary
symmetric functions in the parameter values of the error locations
are Toeplitz, and in particular symmetric. This process is in many
ways reminiscent of describing complete quadrics through various
blowing-ups of the space of
usual quadrics, which is again a space of symmetric square matrices.
Hence, in order to generalize to arbitrary curves, one could
describe some kind of analogy to complete quadrics for curves of
positive genus. Such a generalization, in terms of an object
obtained through various blowing-ups of the secant strata of $C$
inside $\PP$, has been given in Bertram \cite{B}, using vector bundle
language. It is possible that understanding this or similar
objects can give new insight into decoding of Goppa codes.}
\end{rem}

% Section 2
\section{\bf Definitions and basic facts about Goppa codes}
\label{DefBas}

A $q$-ary code of length  $n$  is a subset of the vector space
$F_q^n$, where $F_q$ is a finite field with $q$ elements. A linear
code is a linear subspace of $F_q^n$.  Let  $\overline{F}_q$ be an algebraic
closure of $F_q$. Let  $C$ be a curve of genus $g$ defined over
these fields.
Let  $D$ and $G$ be divisors on $C$ defined over $F_q$, such that
their supports are disjoint. Moreover the support of $D$ consists
of $n$ distinct points $P_1,....,P_n$ of degree one, where
$n=deg(D)$. For a divisor  $M$ defined over $F_q$, denote by
$L(M)$ the set of the zero element and those elements  $f$  of the
function field $F_q(C)$, such that  $(f) + M \ge 0$.  Denote by
$l(M)$ the dimension of $L(M)$ as a vector space over $F_q$.
This dimension is the same as the one obtained if we work over
$\overline{F}_q$.

      Following for example van Lint et al \cite{vLvdG}  we denote by  $C(D,G)$
the code, which is the image of  $L(G)$ in $F_q^n$ under the map:
   $$\phi : f \to\  (f(P_1),.....,f(P_n)).$$
In Pellikaan et al \cite {PSvW}, such a code is called  a WAG, which
is short for
weakly algebraic- geometric code. Moreover one shows there that all
linear codes are WAG. By the theorem of Riemann-Roch  (which remains valid
over finite fields) we have:
      $$dim (C(D,G)) = l(G) - dim(ker f) = l(G) - l(G - D) = $$
%\lineskip
      $$ deg(G) + 1 - g + l(K-G) - l(G-D). $$
As usual $K$ denotes a canonical divisor. Set $m=deg(G)$.  A WAG is
called a SAG (strongly algebraic-geometric code) if the following
composite condition is fullfilled:  $2g-2< m < n$. For a
SAG we observe: $l(K-G) = l(G-D) = 0$, and hence:
                  dim $C(D,G) = m + 1 - g$.
By a generator matrix for a code one means a  $(k \times n)$-matrix,
where the rows constitute a base for the code as a linear space
over $F_q$. A  generator matrix for $C(D,G)$ is:
        $$M  =
\left[\ \vcenter{\halign{$\strut#$\hfill&&\quad\hfill$#$\hfill\cr
f_1(P_1)&f_1(P_2)&\cdots &f_1(P_n)\cr
\cdots \cr
\cdots \cr
f_k(P_1)&f_k(P_2)&\cdots &f_k(P_n)\cr}}\ \right],
$$
where $k=m+1-g$, and ${f_1,....,f_k}$ is a basis for  $L(G)$, both
over $F_q$ and over $\overline{F}_q$.
By a parity check matrix for a code one means a $((n-k) \times n)$-
coefficient matrix of a set of equations cutting out the code as a
subspace of  $F_q^n$ and of $\overline{F}_q^n$.
One denotes by  $C^*(D,G)$ the linear code having the matrix $M$
above as parity check matrix.
Hence $C^*(D,G)$ is the orthogonal complement of $C(D,G)$ and vice
versa. One also says that the codes are dual to each other.
For a WAG defined as $C(D,G)$ consider the  exact sequence of sheaves on C:
      $$0  \to\  \O(G-D) \to\ \O(G) \to\ \O(G)/\O(G-D) \to\  0.$$
The long exact cohomology sequence gives:
$$0  \to\ L(G-D) \to\ L(G) \to\ F_q^n \to\ $$
%\lineskip
$$H^1(C,\O(G-D)) \to\ H^1(C,\O(G)) \to\  0.$$
Moreover by Serre duality :  $H^1(C,\O(G-D))$ is dual to
$L(K+D-G)$, and $H^1(C,\O(G))$ is dual to $L(K-G)$.
For a SAG the long exact sequence reduces to:
\begin{equation} \label{eq:dual}
  0\to\ L(G) \to\ F_q^n \to\ L(K+D-G)^* \to\ 0.
\end{equation}
Here $K$ can (by Riemann-Roch) be chosen such that  $G^* = K+D-G$ has
support disjoint from $D$. We can identify the elements of
$L(K+D-G)$
with differential forms being zero of order the same as order$(G)$
at the points of the support of $G$, and having at most
simple poles at the points of the support of $D$, and no poles
elsewhere. Hence we see that evaluating elements of $L(K+D-G)$ can
be interpreted as evaluating residues of the differential forms
described (after multiplying each value $f(P_i)$ by  a
non-zero value $Res_{P_i}(\eta )$, where $K=(\eta )$  ). Moreover
we have for $f$ in  $L(G)$ and $\omega $ a differential
form as described:
\begin{equation} \label{eq:res}
0  = \Sigma _i Res_{P_i}(f\omega ) = \Sigma _i f(P_i) Res_{P_i}(\omega ).  
\end{equation}

From equations (\ref{eq:dual})  and (\ref{eq:res})  we easily 
conclude that  $C(D,G^*)$ is code equivalent to $C^*(D,G)$.
In the original work by Goppa the code obtained from the divisors
$D$ and $G$ was $C^*(D,G)$, and it was obtained by means of
residues of differential forms. By a Goppa code we will here
simply mean a SAG.  One easily verifies that $C(D,G^*)$ is a SAG if
and only if $C(D,G)$ is so.    For a SAG we see that if $ \{ h_1,
\cdots ,h_{k^*} \} $ is a basis for $L(G^*)$, then a parity check matrix for
$C(D,G)$  is (essentially, after a trivial equivalence operation):
         $$M^*  =
\left[\ \vcenter{\halign{$\strut#$\hfill&&\quad\hfill$#$\hfill\cr
h_1(P_1)&h_1(P_2)&\cdots &h_1(P_n)\cr
\cdots \cr
\cdots \cr
h_{k^*}(P_1)&h_{k^*}(P_2)&\cdots &h_{k^*}(P_n)\cr}}\ \right],
$$  where  $k^* = n-k = n-m+g-1$.
\begin{rem} \label{columns}
{\rm We see that the columns of $M^*$ represent the points of the
support of $D$ if $C$ is embedded into $\PP = \PP^{n-m+g-2}=
\PP (H^0(C,K+D-G)^*)$ by means of sections of $L(G^*)$.}
\end{rem}
\medskip

For $\bf{w}_1$ and $\bf{w}_2$ in $F_q^n$,  let the Hamming
distance $d(\bf{w}_1,\bf{w}_2)$ denote the number of coordinate
positions, in which $\bf{w}_1$ and $\bf{w}_2$ differ; it is
clearly a metric on  $F_q^n$. Let the minimum distance of a code be
the minimum Hamming distance for any pair of codewords. For a linear
code this is easily seen to be the minimum number of non-zero
coordinates (minimum weight) for any non-zero codeword.
    Denote by  $d_1$ the minimum distance of the code $C(D,G)$. If
a code word has weight $d_1$, then there is a divisor $D_1$ with
$D_1 \le D$, and deg($D_1)= n-d_1$, such that  $L(G-D_1) \ne 0$,
Hence  $m-(n-d_1) \ge 0$, that is:  $d_1 \ge n-m$. We denote by
$d$ the integer $n-m$, which is positive since $C(D,G)$ is a SAG.
We call  $d$  the designed minimum distance. We also see that
$C$ is embedded into $\PP (H^0(C,K+H)^*)$, where $d=deg(H)$
(and $H=D-G$). Denote by  $t$  the integer $[(d-1)/2]$. We also call  $t$  the
designed error correcting capacity.  Recall the basic fact:

\begin{rem}
{\rm Let  $N$ be the parity  check matrix of a code. The minimum
distance of the code is equal to $s$ if all choices of  $s-1$
columns of $N$ are independent, and some choice of $s$ columns of
$N$ are dependent.}
\end{rem}
\medskip

Let  $\bf{x}$  be an element (codeword) of $C(D,G) \subset  F_q^n$,
and assume that  $\bf{x}$ is transmitted, and $\bf{y} = \bf{x} +
\bf{e}$  is received.
The difference  $\bf{e}$  is called the error vector.
Denote by  $S(\bf{y})$ the matrix product  $M^*\bf{y}$. Clearly
$S(\bf{y})$ is a vector in $F_q^{k^*}$, and $S(\bf{y})$ = $S(\bf{e})$
 $=\bf{0}$ if and only if  $\bf{y}$ is itself a codeword. $S(\bf{y})$ is
called the syndrome vector of  $\bf{y}$. We can also interpret
$S(\bf{y})$ as a point of  $\PP = \PP^{k^*-1} = \PP^{d+g-2}$.
For each integer  $a$, let the (Hamming) $a$-ball centered at
$\bf{x}$  be the set  of those $\bf{y}$, such that $d(\bf{x},
\bf{y}) \le$ $a$.
The following is immediate from the triangle equality:

\begin{rem}
The restriction of the map $S$:   $F_q^n \to\ F_q^{k^*}$  to any
$t$-ball is injective.
\end{rem}
\medskip
\centerline {\bf{Secant varieties} }
\medskip

Now we view $C$ as any curve defined over the algebraic closure
$\overline{F}_q$, and let $C$ be embedded in some projective space
over this field. Let  $A$ be an effective divisor on $C$, possibly
with repeated points. Let $C_j$ be the $j$'th symmetric product of
$C$, for $j$=1,2,... .

\begin{defn} \label{span}
\begin{itemize}
\item [(a)] We denote by $Span(A)$ the intersection of all hyperplanes
$\H$, such that  we have: $\Sigma _iI(Q_i,C \cap\ \H)Q_i \ge A$  (Here
$I(Q,V_1 \cap\ V_2)$ denotes the usual Bezout intersection number
of two varieties of complementary dimension at a point $Q$).
\item [(b)]We say that $C$ is $k$-spanned if $dim(Span(A))=j-1$,
for all $A$ with $deg(A)=j$, and $j \le k+1$.
\item [(c)]  We set $Sec_j(C) = \cup Span(A)$, where the union is taken
over all  $A$ in $C_j$.
\item [(d)] For a point  $P$ in projective space we set   $h(P) = h$
if  $P$ is contained in  $Sec_h(C) - Sec_{h-1}(C)$.
\end{itemize}
\end{defn}
\begin{prop} \label{emb}
Let $C$ be the curve treated in Section 2, defined over
$F_q$ and embedded into $\PP^{d+g-2}$  by the linear system
$K+D-G$ as described. Then we have:

\begin{itemize}
\item [(a)] $C$ is $(d-2)$-spanned. In particular $C$ is smoothly
embedded if $d \ge 3$.
\item [(b)] If $h(P) = h \le [(d-1)/2] = t$, then there is a unique
effective divisor $A$ with degree at most  $h$, such that  $P$ is
contained in  $Span(A)$.
\end{itemize}
\end{prop}
\begin{proof}
(a) An easy application of Riemann-Roch. Let $A$
be a divisor of degree $j \le d-1$. Set $H = D-G$. Then  $$l(K+H-A)
= 2g-2+d-j+1-g +l(A-H) = d+g-1-j = l(K+H) - j,$$ so $A$ imposes $j$
independent conditions on the linear system.

(b) Assume $P$ is contained in $Span(A_1) \cap\ Span(A_2)$.
If the supports of the divisors $A_1$ and $A_2$ are disjoint,
then we have: $Span(A_1+ A_2)$ = the linear span of  $Span(A_1)
\cup\ Span(A_2)$, so $$dim(Span(A_1+A_2))=1+dim(Span(A_1))+
dim(Span(A_2))-$$

$$ dim((Span(A_1) \cap\  Span(A_2)) \le 1+2(h-1)-1=2h-2.$$
Hence $C$ is not $(2h-1)$-spanned,
and thus not $(d-2)$-spanned, since $h \le [(d-1)/2]$. We leave it
to the reader to modify the argument if the supports of the two
divisors are not disjoint.
\end{proof}

By abuse of notation (See Remark \ref{columns} above) we denote by  $P_i$
column nr. $i$ of the parity matrix $M^*$. Assume that a codeword
$\bf{x}$ is transmitted, $\bf{y}$ is received, and that the error
vector  $\bf{e}$ has weight $h$  with coordinates  $e_1,\cdots ,e_h$
in positions $i_1,\cdots ,i_h$ respectively. We have: The syndrome
$S(\bf{y})$ = $S(\bf{e})$ = $e_1P_1+ \cdots + e_hP_h$. Interpreting
$S(\bf{y})$ as a point of $\PP$, we
then see that  $S(\bf{y})$ is contained in $Sec_h(C)$, and that
$h(S(\bf{y})$)$\le h$.  Moreover it is clear that if  $h \le t$, then
$h(S(\bf{y})$)$ = h$, and that  the "error divisor"
$P_1+\cdots +P_h$ is the unique divisor $A$ of degree at most
$h$ over $\overline F_q$, such that $S(\bf{y})$ is contained in $Span(A)$.
So, error location amounts to finding such a divisor  $A$, given
the point $S(\bf{y})$. A priori we know that this divisor  consists
of distinct points, all of degree 1 defined
over $F_q$, and that even the errors $e_1,\cdots ,e_h$ are in $F_q$.

% Section 3
\section{\bf Vector bundles of rank two on $C$}
\label{Vec2}

We continue using the notation from Section 2.   The following
exposition is to a great extent taken from Lange and Narasimhan
\cite{LN} and Bertram \cite{B}.
Let $Ext_{\O_C}(H,\O_C)$ be the set of isomorphism classes of
exact sequences  $$(e): 0 \to\ \O_C \to\  E \to\ H \to 0.$$
The map $\O_C \to\  E$ is denoted by $f$ and the map $E \to\ H$
by $g$. The zero element $(e_0)$ corresponds to the case of a split
exact sequence.
Here $\O_C$ as usual denotes  the structure sheaf on $C$, and $H$
is the fixed line bundle or invertible sheaf $D-G$, see
Section 2 (by abuse of notation we do not distinguish between
the divisor $H = D-G$, or the invertible sheaf or line bundle, of
which the divisor corresponds to a global section). The middle
term  $E$ is a locally free sheaf, or vector bundle, of rank 2.
Standard cohomology theory and Serre duality give:
$$Ext_{\O_C}(H,\O_C) = Ext_{\O_C}(\O_C,-H) = H^1(C,-H) = H^0(C, K+H)^*.$$
Hence  $\PP(Ext_{\O_C}(H,\O_C)) = \PP(H^0(C, K+H)^*) = \PP$.
This means that (up to isomorphism and a multiplicative factor) the
points of our well-known projective space $\PP$ described in
Section 2 are identified with extensions as described.

\begin{defn} \label{s-inv}
%$\bf{Definition}$ $\bf{(3.1)}$
Let $E$ be a rank two vector bundle on $C$.
  \begin{itemize}
\item [(a)] Denote by  $s(E)$ the integer
$$deg(E)-2max(deg(L))=2min(deg(M))- deg(E),$$
where the maximum is taken over all line subbundles $L$ of $E$ and
the minimum is taken over the quotient line bundles $M$ of $E$.
\item [(b)]  $E$ is called stable if $s(E) > 0$; semistable if
$s(E) \ge 0$;  unstable if $s(E)< 0$.
\item [(c)] For an extension $(e)$ as above we set $s((e))=s(E)$, where $E$
is the middle term. An extension is  called stable (semistable,
unstable) if the middle term E is so.
 \end{itemize}
\end{defn}
The definitions of stable and semistable coincide if
$deg(E)$ is odd.
\medskip
For the zero element  $(e_o)$ we observe:   $s(E) = -d$, and for
all non-split extension $(e)$ we have $d \ge s((e)) \ge 2-d$.
Moreover, if $M$ is a  quotient line bundle of $E$ of minimal
degree $(s+d)/2 \ge 1$, with quotient map $h$, then the composition
of $h$ and $f$  is non-zero:
                  $$ \O_C \to\ E \to\ M.$$
Hence $M$ is isomorphic to $\O_C(A)$, for an effective divisor $A$
of degree $(s+d)/2$. This implies again that (identifying $(e)$
with its corresponding point of $\PP$) the point $(e)$ is contained
in the kernel of the map:
$$ Ext^1_{\O_C}(H,\O_C) \to\  Ext^1_{\O_C}(H,\O_C(A)),$$
that is, in the kernel of:
$$ H^0(C, K+H)^* \to\ H^0(C, K+H-A)^*.$$
This observation has many consequences. First we see that the
set of points in $\PP$ with $s$-value $2-d$ are precisely those
that  represent  bundles with a quotient bundle of type $\O_C(Q)$,
for some point $Q$ on $C$. If we assume that $d \ge 3$, then $C$
is smoothly embedded by Proposition \ref{emb}. Then we can identify
$C$ with its embedded image in $\PP$, and the point $Q$ on C then
corresponds to a bundle extension with a quotient bundle isomorphic
to $\O_C(Q)$. Moreover it is clear that the observation above is
equivalent to:
$(e)$ is contained in $Span(A)$.

\begin{rem} \label{spana}
{\rm Arguing in a dual way, we get that if (e) is contained in
Span(A), then the line bundle corresponding to $H-A$
is a subbundle of E.}
\end{rem}
Summing up, we now formulate the following result,
which is practically identical to Proposition 1.1. of Lange and
Narasimhan \cite{LN}
(Recall the functions  $h$  and  $s$  from $\PP$  to $\bf{Z}$,
introduced in Definition \ref{span}, d.) and Definition \ref{s-inv}, (a),
respectively).

\begin{prop} \label{s and h}
Let  $P$ be a point of $\PP$. Then   $s(P) = 2h(P) - d$.
In particular  $P$ is a unstable point  (semistable, stable) if and
only if   $h(P)<d/2$  $ (h(P) \ge d/2,$ $ h(P) > d/2)$.
\end{prop}
We are now able to formulate the main result of the paper:
\begin{thm} \label{mainthm}
Let  $P = S(\bf{y})$ be the syndrome of a received message using the
code $C(D,G)$. Then $P$ is the syndrome of some error vector with
weight at most the designed error correcting capasity $t = [(d-1)/2]$
only if $P$ is an unstable point.
Moreover in that case the process of error location is
reduced to finding the error divisor $A$ among global sections of
the unique quotient line bundle of degree $h$ of
the vector bundle $E(P)$ of rank two, appearing as the middle term in the
extension corresponding to $P$.
\end{thm}
\begin{proof}
This follows directly from Proposition \ref{s and h} and the argument above.
\end{proof}

\begin{rem} \label{iff}
{\rm The ``only if'' in the theorem can be replaced by ``if and only if'' 
if we  define the syndrome map over $\overline{F}_q$.}
\end{rem}
\begin{rem} \label{GIT}
{\rm The definitions of stable, semistable, unstable can be viewed as
special cases of more general definitions of these concepts in the
setting of Geometric Invariant Theory (GIT), which again is an
essential tool in building moduli spaces parametrizing various
objects. The spaces arise as quotients of various group actions.
In order to get quotients with good properties one usually has
to disregard certain "bad objects", which are the unstable ones.
In our case the relevant construction is that of $\mathcal{SU}$$_2(C)$,
the moduli space of isomorphism classes of vector bundles of rank
two on $C$  (See Gieseker \cite{Gi}, p. 51-52). This has dimension
$4g-3 = g + (3g-3)$ , where the sum decomposition corresponds to
$g=dim(Jac(C))$ degrees of freedom to choose a line bundle $H$,
and $3g-3= dim (\mathcal{SU}$$_2(C,H))$  degrees of freedom to choose
the rank two bundle with determinant $H$ (where $\mathcal{SU}$$_2(C,H)$
 = modulo space of rank two bundles with determinant $H$). Moreover,
for $\mathcal{SU}$$_2(C,H)$ there are essentially only 2 cases;
$deg(H)$ odd and $deg(H)$ even, since tensoring a rank two bundle
with a line bundle gives rise to an isomorphism between two such
spaces with determinants $H$ with degrees of equal parity.  One can
show that the natural map:
             $$\PP - Sec_t(C) \to\ \mathcal{SU}_2(C,H)$$
is  birational if  $d = deg(H) = 2g-1$, and that it maps
birationally on to the $\theta$-divisor if $d = deg(H) = 2g-2$,
where  the $\theta$-divisor parametrizes the rank two-bundles
with a global section.

   One observes that the situation in coding theory in a certain
way is complementary to that of applying GIT to build
$\mathcal{SU}$$_2(C,H)$, since the good(syndrome) points in coding
theory are the bad ones for GIT, and vice versa. On the other hand
the issue for those who work with moduli spaces is often precisely
what to do with the unstable points, so the focal point of the
theories are still in a certain sense overlapping. One is for
example interested in blowing up various $s$-negative strata of
$\PP$ to obtain compactifications of $\PP - Sec_t(C)$  and
$\mathcal{SU}$$_2(C,H)$  with desired properties. See Bertram \cite{Bdiff}
and \cite{B}. One could hope that insight in such compactifications
could be instrumental in understanding algorithms for decoding of
Goppa codes.
   One can also ask: Is it possible that the minimum distance of
$C(D,G)$ exceeds $d$, or weaker: Given a reasonable small integer
$k_0$, like 2, 3 or 17; is there a positive limit $s_0$, such that
if $s(P) \le s_0 = 2h_0 - d$, then there are at most $k_0$ divisors
$A$ (even over $\overline{F}_q$) of  degree $h_0$, such that  P is
contained in $Span(A)$? For practical purposes this would in some situations
be almost as good as unique decoding. In vector bundle language one
is then interested maximal sublinebundles of $P$ corresponding to
divisors of type $H-A$. The question of such maximal subbundles is
the main issue in Lange and Narasimhan \cite{LN}.}
\end{rem}
\begin{example} \label{ratcurv}
{\rm Assume $g=0$. Then $C$ is mapped into $\PP = \PP^{d-2}$ as a curve
of degree $d-2$, that is as a rational normal curve. It is
well-known that on $C = \PP^{1}$ all rank two bundles of degree $d$
split as $\O(a) \oplus\ \O(b)$, with $a+b = d$. If $d$ is odd, we
then see the largest possible $s$-value is $-1 = d- 2[(d+1)/2]$.
Hence all rank two-bundles are unstable, corresponding to the fact
that $\PP = Sec_t(C)$ over $F_q$.}
\end{example}

\centerline {\bf Acknowledgements}
\medskip
I am grateful to Emma Previato for informing me about the manuscript 
\cite{BC}, and to the authors of \cite{BC} for correcting a mistake
in an old version of my paper. 
\medskip

%\begin{thebibliography}{[1]}
\centerline{\bf References}
%\begin{thebibliography}{[E-L-M-S]}
\begin{itemize}
\bibitem[1]{A} M.~Atiyah,  Complex fibre bundles and ruled
    surfaces, \textit{Proc. London Math. Soc.}, \textbf{5} (1955), 407-434.

\bibitem[2]{Bdiff} A.~Bertram, Moduli of rank-2 vector 
bundles, theta divisors, and the geometry of curves in projective
space, \textit{Journal of Diff. Geometry }, \textbf{35} (1992), 429-469.

\bibitem[3]{B} A.~Bertram, Complete Extensions and their Map
to Moduli Space,  \textit{London Math. Soc. Lect. Notes Series}, \textbf{179}  
(1992), 81-91.

\bibitem[4]{BC} T.~Bouganis, D. Coles, A Geometric View of Decoding AG
  Codes,  Manuscript, 10 pages  (2002).

\bibitem[5]{D} I.~Duursma, Algebraic decoding using special
divisors, \textit{IEEE Trans. of Info. Theory},  \textbf{39(2)}
(1993), 694-698.

\bibitem[6]{FR} G.L. ~Feng, T.R.N.~Rao, Decoding 
Algebraic-Geometric Codes up to the Designed Minimum Distance, 
\textit{IEEE Trans. of Info. Theory},  \textbf{39(1)} (1993), 37-45.

\bibitem[7]{G1} V.D.~Goppa, Algebraic-Geometric codes, 
\textit{Math. USSR Izvestiya}, \textbf{21} (1983), 75-90.

\bibitem[8]{G2} V.D.~Goppa, \textit{Geometry and Codes}, 
Mathematics and its applications, Soviet series, \textbf{21},  Kluwer
Ac.Publ., (1989).

\bibitem[9]{Gi} D.~Gieseker, Geometric invariant theory and 
applications to moduli problems, \textit{Springer Lect. Notes in Math.},  
\textbf{996} (1983), 45-63.

\bibitem[10]{JLJH} J. ~Justesen, K.J. ~Larsen, H. ~Elbr\"ond Jensen, 
and T. ~H\"oholdt, Fast decoding of codes from algebraic plane 
curves, \textit{IEEE Trans. of Info. Theory},  \textbf{38(1)} (1992), 111-119.

\bibitem[11]{LN} H.~Lange, M.S.~Narasimhan, Maximal subbundles 
of rank 2 vector bundles on curve, \textit{Math.Ann.}, \textbf{266} (1983), 55-72.

\bibitem[12]{P} R.~Pellikaan, On a decoding algorithm for codes on 
maximal curves, \textit{IEEE Trans. of Info. Theory}, \textbf{35(6)} 
(1989), 1228-1232.

\bibitem[13]{PSvW} R. ~Pellikaan, B.Z. ~Shen, and G.J.M. ~van Wee  
Which linear codes are algebraic geometric?, \textit{IEEE Trans. of
Info. Theory},  \textbf{37(3)} (1991), 583-682.

\bibitem[14]{SV} A.N. ~Skorobogatov and S.G. ~Vladut,  On the
    Decoding of Algebraic-Geometric Codes,  \textit{IEEE Trans. of
  Info. Theory } \textbf{36(5)} (1990), 1051-1060.

\bibitem[15]{vLvdG} J. ~van Lint and G. ~van der Geer,  \textit{Introduction to Coding Theory and Algebraic Geometry}, DMW Seminar \textbf{12}, Birkhauser (1995).
\end{itemize}

%\end{thebibliography}

\end{document}